\documentclass[amsmath]{kluwer}    % Specifies the document style.

\newdisplay{guess}{Conjecture}

\begin{document}                                                                                   
\begin{article}
\begin{opening}         
\title{Smectic Liquid Crystals: Materials with One-Dimensional, Periodic Order} 
\author{Randall D. \surname{Kamien}}  
\author{Christian D. \surname{Santangelo}}
\runningauthor{Randall D. Kamien and Christian D. Santangelo}
\runningtitle{Smectic Liquid Crystals}
\institute{Department of Physics and Astronomy, 209 South 33$^{\sl rd}$ Street, University of Pennsylvania, Philadelphia, PA 19104-6396, USA}
\date{6 December 2005; revised 16 January 2006}

\begin{abstract}
Smectic liquid crystals are materials formed by stacking deformable, fluid layers.  Though smectics prefer to have flat, uniformly-spaced layers, boundary conditions can impose curvature on the layers.  Since the layer spacing and curvature are intertwined, the problem of finding minimal configurations for the layers becomes highly nontrivial.  We discuss various topological and geometrical aspects of these materials and present recent progress on finding some exact layer configurations.  We also exhibit connections to the study of certain embedded minimal surfaces and briefly summarize some important open problems.
\end{abstract}

\end{opening}

\section{Liquid Crystals: Phases in Between Liquid and Solid}
First characterized at the end of the 19th century, liquid crystals have gone from chemical curiosities to important optical and electronic devices~\cite{degennes}.  In fact, this article relied on liquid crystal technology to be produced!  
As their odd name implies, liquid crystals have some properties of a crystal and some of a liquid.  Whether {\sl nematic}, {\sl cholesteric}, or {\sl smectic}, all these phases flow (like a fluid), at least in some directions, while transmitting torques and forces (like a crystal) in others.  The phase we consider
here is the {\sl smectic-A} phase.  A crystal is characterized as a periodic structure (the simplest of which are lattices).  A $d$-dimensional crystal is a periodic structure in $d$ dimensions with
$d$ independent periodicities.  However, it is possible to have a lower number of periods than
than number of dimensions.  Focussing on $d=3$ dimensions, a one-dimensional periodic structure
is a {\sl smectic}, a two-dimensional periodic structure is a {\sl columnar phase}, a three-dimensional
structure is a crystal, and a phase with no periodicity is a fluid.

In this article we will consider smectic phases and discuss minima of the free energy
\begin{equation}\label{functional}
F = \frac{1}{2}\int_V d^3\!x\,\left\{B\left(\vert\nabla\Phi\vert -1\right)^2 + K_1\left(\nabla\cdot{\bf N}\right)^2\right\}
\end{equation}
where the density of material is $\rho({\bf x})\propto \cos\left[2\pi\Phi({\bf x})/a\right]$, $a$ is the one-dimensional spacing of the layers and ${\bf N}=\nabla\Phi/\vert\nabla\Phi\vert$ is a unit vector field, normal to the level sets of $\Phi$.  The constants $B$ and $K_1$ are the bulk and bend modulus, respectively, and the integral is over the sample volume $V$.  In particular we are interested in minima with non-trivial topology on domains which are not simply-connected.  Before discussing those solutions, we will review the origins of the functional in (\ref{functional}).

A one-dimensional periodic structure in three dimensions can be represented as an $L_2$ summable density
wave in $\mathbb{R}^3$
\begin{equation}
\rho({\bf x}) = \sum_{j=-\infty}^\infty \rho_j e^{ij{\bf q}\cdot {\bf x}}
\end{equation}
where $\rho_j^*=\rho_{-j}$, so that $\rho({\bf x})$ real, and ${\bf q}$ is a vector pointing along the direction
of the modulation with magnitude $2\pi/a$, where $a$ is the spacing between the ``layers''.    Each Fourier mode $\rho_j$ for $j\ne 0$ is an {\sl order parameter} for the smectic: if $\rho_j =0$ for all
$j\ne 0$, then the material is a fluid, while if some $\rho_j$ is nonvanishing, it is a smectic.  Since $a$ is the spacing between the layers, it is sufficient to keep track of just $\rho_1$ since the presence of the fundamental mode is enough to signal the order.  Truncating (1) to the first non-trivial mode and writing $\rho_1=\vert\rho_1\vert e^{-2\pi i u/a}$~\cite{chaikin} we have
\begin{equation}
\rho({\bf x}) = \rho_0 + 2\vert\rho_1\vert \cos\left[{2\pi\over a}\left(\hat q\cdot{\bf x} -u\right)\right]
\end{equation}
To study deformations and fluctuations of this periodic structure, we promote $u$ from a constant
to a scalar field $u({\bf x})$.  This results in an {\sl Eulerian} description of the material, {\sl i.e.} the
layers sit at surfaces of constant $\Phi({\bf x})\equiv \hat q\cdot{\bf x} - u({\bf x})$ or, equivalently, at the peaks of the density wave.  An alternative description is {\sl Lagrangian}, in which we label each fiducial plane defined by $\hat q\cdot {\bf x}_n =n$ and
consider deformations about these.  When the deformations are small they can be represented
as a height function or graph over each plane.  If we choose our coordinates so that $\hat q=\hat z$, then
\begin{equation}
n= h_n(x,y) - u\left[x,y,h_n(x,y)\right]
\end{equation}
However, a unique solution for $h_n$  requires that the deformed surface
not have any overhangs and, moreover, awkwardly mixes spatial coordinates with the Lagrangian
coordinate $n$ which labels the original surface.  Despite this, as we will see subsequently, in restricted cases the Lagrangian representation is a useful tool.

There are three essential actors in the free energy of smectics: the layer compression or strain, the
bending or mean curvature, and the topology of the surfaces characterized by the intrinsic or Gaussian curvature.  
Using the phase field $\Phi$, related to the density $\rho$ via
$\rho\propto\cos\left(2\pi\Phi/a\right)$, the unit layer normal
is 
\begin{equation}
{\bf N} = \frac{\nabla\Phi}{\vert\nabla\Phi\vert}
\end{equation}
At a point on a surface, the two principle radii of curvature $R_1$ and $R_2$ determine the
two curvatures $\kappa_i\equiv 1/R_i$.  These two curvatures can be combined into the
mean curvature $H={1\over 2}\left(\kappa_1+\kappa_2\right)$ and the Gaussian curvature
$K=\kappa_1\kappa_2$.  The former is a measure of the bending of the layers in space, while
the latter is an intrinsic quantity which is a local measure of geometric distortion and can be
applied globally to gather topological information.   Fortunately, these quantities may
be expressed in terms of the layer normal:~\cite{kamien2}
\begin{eqnarray}
H & =& {1\over 2}\nabla\cdot{\bf N}\bigg\vert_{\Phi(x,y,z)=c}\\
K & = &{1\over 2}\nabla \cdot \left[   {\bf N} \left( \nabla\cdot {\bf N} \right) -\left( {\bf N} \cdot \nabla \right) {\bf N}\right]\bigg\vert_{\Phi(x,y,z)=c}
\end{eqnarray}

\noindent where the expressions are evaluated on one of the level sets of $\Phi$.  By working in this set of 
coordinates, it is possible to free these geometric quantities from the surfaces and, as we will
show, allows us to relate different terms in the free energy.  

Because smectics have a preferred layer spacing, $a$, there is necessarily a strain energy associated
with changing this spacing.  While the two curvatures are natural from the geometric point of view, an expression for the compression strain is not as well defined.  Returning to the family of equally spaced layers, we have $na=\Phi({\bf x}_n +an {\bf N})$ and so, differentiating with respect to $n$, find that
\begin{equation}
1={\bf N}\cdot\nabla\Phi({\bf x}_n + an {\bf N})
\end{equation}
which implies that $\vert\nabla\Phi\vert=1$.  Thus the strain, when expressed in terms of $\Phi$ must vanish when this condition holds~\cite{klemanbook}.  For instance, we can define the strain $u_{zz}$ as
\begin{equation}
u_{zz}  \equiv 1-\left\vert\nabla\Phi\right\vert  = 1-{\bf N}\cdot\nabla\Phi
\end{equation}
or an alternate strain $\tilde u_{zz}$ as
\begin{equation}
\tilde u_{zz} \equiv {1\over 2}\left[1-\left(\nabla\Phi\right)^2\right]
\end{equation}
both of which vanish when $\vert\nabla\Phi\vert=1$.  Moreover, 
in terms of $u$, the Eulerian displacement field, it follows that $u_{zz} \approx \partial_z u + {\cal O}\left(u^2\right)$
and $\tilde u_{zz} \approx \partial_z u + {\cal O}\left(u^2\right)$ are identical to linear order in $u$ and differ only
in their nonlinear terms.  Therefore, when considering the quadratic compression energy:
\begin{equation}
F_{\rm c} = {1\over 2}\int d^3\!x\, B u_{zz}^2
\end{equation}
the harmonic term is independent of the specific choice of $u_{zz}$ -- the linear term will
always be $\partial_z u$.  While $\tilde u_{zz}$ arises naturally from a Landau theory description
of the nematic  to smectic transition~\cite{degennes72}, there is no compelling reason to choose it over
$u_{zz}$ -- both are invariant under rotations and vanish for equally spaced layers.  The constant $B$ is known as the {\sl bulk modulus} and controls compression modes of the layers.  As we will show in the following, the use of $u_{zz}$ allows us to further develop the
geometric theory of smectics \footnote{Indeed $u_{zz}=1- \sqrt{1-2\tilde u_{zz}}= \tilde u_{zz} + {\cal O}\left(\tilde u_{zz}^2\right)$ and so energies based on $u_{zz}$ will differ from those with $\tilde u_{zz}$ {\sl only} by anharmonic terms. }.  In a similar vein, one could argue that the mean curvature $H$ could be replaced by ${1\over 2}\nabla^2\Phi$.  However, the nonlinear form $H={1\over 2}\nabla\cdot{\bf N}$ {\sl unambiguously} has the Laplacian as the first term in an expansion in powers of $u$, while there are many nonlinear generalizations of the Laplacian, including the full mean curvature.

Putting this all together we have~\cite{santangelo2}
\begin{equation}\label{fullfreeenergy}
F ={B\over 2}\int d^3\!x\,\left\{\left[1-{\bf N}\cdot\nabla\Phi\right]^2 + \lambda^2\left(\nabla\cdot{\bf N}\right)^2 \right\},
\end{equation}
where the last term is the contribution of the layer curvature to the energy, and is quadratic in the radii of curvature.  The constant $\lambda=\sqrt{K_1/B}$ has units of length and sets the natural lengthscale for deformations.  For the sake of consistency, one should add to this expression a term of the form
\begin{equation}\label{gauss}
F_K = \bar{K} \int_V d^3\!x\,\nabla \cdot \left[ (\nabla \cdot {\bf N}) {\bf N} - \left( {\bf N} \cdot \nabla \right) {\bf N} \right],
\end{equation}
coming from the Gaussian curvature.  Since this term is a total derivative, it can be reduced to an integral on the boundary and, therefore, cannot play a role in the minimization of the energy that determines the configuration of the layers.  The usual ground state which minimizes (\ref{fullfreeenergy}) depends on an arbitrary unit vector $\bf n$: $\Phi= {\bf n}\cdot{\bf x}$, corresponding to equally spaced layers along 
the $\bf n$ direction.  The free energy of this solution vanishes since ${\bf N}={\bf n}$ is constant and
$\vert\nabla\Phi\vert=\vert{\bf n}\vert=1$. 

It will become clear why we have chosen this nonlinear form for $u_{zz}$.  In the following section, however, we will linearize this theory -- a process which results in an unambiguous  form for $u_{zz}$.

Only in very few cases can we satisfy the terms in the energy individually.  The difficulty of this problem has to do with the nontrivial way that the terms in equation~(\ref{fullfreeenergy}) interact.  It is
instructive to consider
a family of evenly spaced layers ${\bf x}_n(s_1,s_2)$ parameterized by $(s_1,s_2)$, with ${\bf x}_n={\bf x}_0 + na{\bf N}$.
At any point $(s_1,s_2)$ on ${\bf x}_0$, we have the two principal curvatures $\kappa_{1,0}$ and $\kappa_{2,0}$.  The radii of curvature of the nth layer are just
$R_{i,n} = R_{i,0} + na$
from which it follows that 
\begin{equation}
\kappa_{i,n} = \frac{\kappa_{i,0}}{1+na\kappa_{i,0}}
\end{equation}
From the definitions of $H$ and $K$
it follows that the curvatures of the nth layer are:
\begin{eqnarray}
K_n &=& \frac{K_0}{1+2naH_0 + n^2a^2K_0}\\
H_n &=& \frac{H_0 + naK_0}{1 + 2naH_0 + n^2a^2 K_0}
\end{eqnarray}
Thus if the layers are evenly spaced, the only way for the bending to vanish for all layers ($H_n=0$) is for $K_0$ (and thus $K_n$) to vanish~\cite{didonna}.  This interplay between even spacing, the Gaussian curvature, and the mean curvature is a recurring motif in this work and will be discussed further in the following.
It is clear that in order to have configurations with non-vanishing curvature, the layer spacing
must necessarily vary.

\subsection{Topological Defects}

The ground state or equilibrium configuration of the field $u$ subject to boundary conditions
is the goal of our analysis.  While the boundary conditions at infinity can usually be taken to be
$u=0$ or $\nabla u=0$, we may introduce additional boundary conditions via topological defects~\cite{chaikin}.  In particular, because $\Phi$ is a phase field, it can wind by multiples of $2\pi$ around defects.  In the
so-called ``London limit'', we simply remove a line from $\mathbb{R}^3$ around which the phase winds.  It is possible to instead allow the density $\rho_1$ to vary in space: the solution to the full equation has
$\rho_1$ vanishing at the core of the defect rendering the energy finite.  This extra complication is 
tangential to our discussion and we will use the London theory in the following.  

Much of the physical understanding of defects comes from the linearized form of the theory or, equivalently, the quadratic energy~\cite{klemanbook,degennes}:
\begin{equation}
F_l = \frac{B}{2} \int d^3x~\left[ \left(\partial_z u \right)^2 + \lambda^2 \left( \nabla_\perp^2 u \right)^2 \right]
\end{equation}
The first term is the linearized compression strain and measures the stretching along the layer direction, taken to be $\hat z$.  We define $\nabla_\perp = \hat x\partial_x +\hat y\partial y$, the gradient in the plane perpendicular to $\hat z$, and see that the second term is the Laplacian in that plane, which we
recognize as the linearized form of the mean curvature, $H$.
It is important to note that the original theory was quadratic in the strain $u_{zz}$ and the 
curvature $H$.  There is no reason, physical or otherwise, that
this truncation should be valid -- it is used primarily because of its solvability.  In addition, the ambiguity arising from the freedom to choose different forms of the nonlinear free energy does not arise in the linear theory.  In order for the density to be single-valued, we have
\begin{equation}\label{integral}
\oint_\Gamma \nabla u\cdot d\ell = na
\end{equation}
for some $n\in\mathbb{Z}$ whenever the contour $\Gamma$ goes around the defect.  The combination $-na=b$ is called the {\sl Burgers vector}~\cite{chaikin}.  This can
be promoted to differential form by introducing a defect density ${\bf m}({\bf x})$ which is a vector
field pointing along the defect line~\cite{degennes}.  The defect line is along the curve ${\bf r}(s)\in\mathbb{R}^3$ with local unit tangent  ${\bf t}(s)=d{\bf r}/ds$ and we write
\begin{equation}
{\bf m}({\bf x}) = na\int_0^L \delta^3\left[{\bf x} - {\bf r}(s)\right] {\bf t}(s)
\end{equation}
It follows that $\nabla\cdot {\bf m} =0$ for defect lines that either run off to infinity or are closed loops.
From Stoke's theorem, we have the local condition:
\begin{equation}
\nabla\times\nabla u = {\bf m}
\end{equation}
To find the ground state configuration of the smectic in the presence of defects, we write
$u$ as the sum of a regular part $\sigma$ and a singular part $\hbox{\boldmath$\tau$}$ so that
$\partial_i u = \partial_i\sigma + \tau_i$ with $\nabla\times\hbox{\boldmath$\tau$} ={\bf m}$~\cite{degennes}.  This is most easily done
in Fourier space where
\begin{equation}
\hbox{\boldmath$\tau$}({\bf q}) = {-i {\bf q}\times {\bf m}({\bf q})\over q^2}
\end{equation}
where we have chosen $\nabla\cdot\hbox{\boldmath$\tau$} =0$, absorbing any of the longitudinal parts of $\hbox{\boldmath$\tau$}$ into
$\sigma$.  Having enforced the boundary conditions via the introduction of the singular part of $u$, we insert the expression for $\nabla u$ into the free energy and minimize over the remaining variable, $\sigma$.  This results in
\begin{equation}
0=\frac{\delta F_l}{\delta \sigma}  = B\left\{\left[q_z^2 +  \lambda^2q_\perp^4\right] \sigma - q_z\left(1-\lambda^2q_\perp^2\right) \frac{\hat z\cdot(q\times m)}{q^2}\right\}
\end{equation} 
so that
\begin{equation}\label{tw1}
\sigma(q) = \frac{q_z \left(1-\lambda^2q_\perp^2\right)\left(q_x m_y - q_y m_x\right)}{q^2\left(q_z^2 + \lambda^2q_\perp^4\right)}
\end{equation} 
and 
\begin{equation}
F_l  = {B\lambda^2\over 2}\int \frac{d^3\!q}{(2\pi)^3} {\vert \hat z\cdot(q\times m)\vert^2\over q_z^2+\lambda^2q_\perp^4} 
\end{equation}

Because of the anisotropy of the smectic, we can have defects which are parallel to the special, $\hat z$ axis and defects which are perpendicular to it.  The former are called {\sl screw dislocations} while the 
latter are known as {\sl edge dislocations}~\cite{klemanbook}.  Defects can, of course, have a mixed character when
they are oblique to the $\hat z$ axis.  For an edge dislocation pointing along the $\hat y$ axis,
$m=na\hat y \delta(x)\delta(z)$ or $m(q)=(2\pi) na \hat y\delta(q_y)$. In order to calculate the energy we will need a long length cutoff $L_y/2$ for $q_y$ and a short distance cutoff $q_{\rm max} = \pi/\xi$ for $q_z$.  We have 
\begin{eqnarray}
F_{\rm edge} &=& {B\lambda^2n^2a^2\over 4\pi} \int {dq_x dq_z} {q_x^2\over q_z^2+\lambda^2q_x^4} \delta(q_y=0)\nonumber\\ 
&=& {BL_y\sqrt{\lambda} n^2a^2\over 4\sqrt{2}} \int_{-q_{\rm max}}^{q_{\rm max}} dq_z {1\over\sqrt{\vert q_z\vert}}\nonumber\\
&=& {BL_y n^2a^2\sqrt{\pi}\over \sqrt{2}}\sqrt{{\lambda}\over { \xi} }
\end{eqnarray}
so that the energy per unit length ($F/L_y$) is finite~\cite{degennes,klemanbook}.  For a screw dislocation, $m=na\hat z\delta(x)\delta(y)$ and
the free energy vanishes meaning that the smooth part of $u$ can ``screen'' the strain
from the boundary at the defect.  The spatially dependent solutions for the edge and screw dislocations are
\begin{eqnarray}
\nabla u_{\rm edge} &=&  {na\over 4\sqrt{\pi\lambda \vert z\vert}}\exp\left\{-{x^2\over 4\lambda \vert z\vert}\right\}\left[1,0,-{x\over z}\right]\\ \label{eq:linearscrew}
\nabla u_{\rm screw} &=& {na\over 2\pi \left(x^2+y^2\right)}\left[y,-x,0\right]
\end{eqnarray}
respectively.

\subsection{The Edge Again}
It is instructive to notice that the result for edge dislocation is simply the gradient of 
\begin{equation}
u(y,z) =-\hbox{sgn}(z){b\over 2}{1\over\sqrt{2\pi}}\int_{-\infty}^{x\over \sqrt{2\lambda\vert z\vert}} dt\, e^{-t^2/2}
\end{equation}
where we will now switch to the Burgers vector $b=-na$.  For $z>0$, we have
$\left(\partial_z -  \lambda\partial_x^2\right)u=0$ while for $z<0$, $\left(\partial_z + \lambda\partial_x^2\right)u=0$.   There is a discontinuity between the upper and lower half planes which precisely
enforces the condition (\ref{integral}).  To see this, we take the clockwise path in the $yz$ plane, starting at $(y,z)=(-\infty,0^-)$, going to $(+\infty,0^-)$, to $(+\infty,0^+)$ and back to $(-\infty,0^+)$.  We have 
\begin{eqnarray}
\int_\Gamma \nabla u\cdot d\ell &=& \left[u(-\infty,0^+)-u(\infty,0^+)\right]
+\left[u(\infty,0^+)-u(\infty,0^-)\right] \nonumber \\&&+\quad \left[u(\infty,0^-)-u(-\infty,0^-)\right]\nonumber\\
&=& {b\over 2} -b -{b\over 2} = -b
\end{eqnarray}
Using this result, we can calculate the energy directly by noting that
\begin{eqnarray}
F_l &=& {B\over 2}\int_M d^3\!x\,\left\{\left[\partial_z u \mp \lambda\nabla_\perp^2 u\right]^2 \pm 
2\lambda\partial_z u\nabla_\perp^2u\right\}\\
&=&{B\over 2}\int_M d^3\!x\,\left\{\left[\partial_z u \mp \lambda\nabla_\perp^2 u\right]^2\right\}\nonumber\\ &&\quad\pm 
B\lambda \int_{M} d^3\!x\, \left\{\nabla_\perp\cdot\left(\partial_z u\nabla_\perp u\right) - {1\over 2}\partial_z\left(\nabla_\perp u\right)^2\right\}\nonumber
\end{eqnarray}
We take $M$ to be all of $\mathbb{R}^3$.  In the upper half plane, we pick the minus sign, while in the lower half plane the plus sign so that
the first integral over $M$ vanishes and the second becomes a surface integral over $\partial M$, the $xy$ plane.  We see then that
the energy is
\begin{eqnarray}
F_l &=& {B\lambda} L_y {b^2\over 16\pi\lambda\xi}\int_{-\infty}^\infty dx\, \exp\left\{{-x^2\over 2\lambda\xi'}\right\}\nonumber\\
&=& {B L_y b^2\over 8\sqrt{2\pi}}\sqrt{\lambda\over\xi'}
\end{eqnarray}
where, again, we have introduced a short-distance $z$ cutoff, $\xi'$.  By choosing $\xi'=\xi/(64\pi^2)$
the two expressions for the energy can be made to agree.  The important point is that the energy
has the same dependence on the physical parameters, $b$, $B$ and $\lambda$.  

The insight of Bogomoln'yi, Prasad, and Sommerfield, was to turn this calculation around: by rewriting the free energy as a perfect square plus a total derivative we can bound the energy from
below simply through the boundary conditions~\cite{BPS1,BPS2}.  Moreover, when there are solutions of the
lower order differential equation in the perfect square satisfying the appropriate boundary conditions, the
energy arises {\sl solely} from the boundary conditions.  Such states are called BPS states and their energy is sometimes called ``topological'' because of its sole dependence on the boundary conditions.

In the next section we consider the rotationally invariant model with this decomposition in mind.

\section{Rotational Invariance: From Here to Nonlinearity}
The quadratic approximation employed in the previous section has the drawback that it does not preserve the invariance of the free energy under rigid rotations of the layered structure.  This rotation invariance is preserved by the higher order terms that would appear in an explicit expansion of the energy in the deformation $u(x,y,z)$.  As it happens, these nonlinear terms are crucial for understanding a number of physical phenomena not captured by the quadratic approximation, yet they render the energy too complex for a straightforward analytical solution.

\subsection{BPS Solutions}
In the full free energy, there are also exact minima that satisfy a BPS flow equation~\cite{santangelo1,santangelo2}.   We start
with the nonlinear theory (\ref{fullfreeenergy}):
\begin{eqnarray}
F ={B\over 2}\int d^3\!x\,\left\{\left[1-{\bf N}\cdot\nabla\Phi\right]^2 + \lambda^2\left(\nabla\cdot{\bf N}\right)^2\right\}\nonumber
\end{eqnarray}
Following the discussion in the last section, we complete the square yielding:
\begin{equation}\label{withcrossterm}
F = \frac{B}{2} \int d^3x~\left[ \left(1-| \nabla \Phi| - \lambda \nabla \cdot {\bf N} \right)^2 + 2 \lambda (1-|\nabla \Phi|) \nabla \cdot {\bf N} \right]
\end{equation}
While the first cross term is a total derivative ($\nabla\cdot{\bf N}$), the second is not.  However, we
write 
\begin{eqnarray}
\vert\nabla\Phi\vert\nabla\cdot{\bf N} &=&
\left[{\bf N}\left(\nabla\cdot{\bf N}\right)\right]\cdot\nabla\Phi\nonumber\\
&=&\left[{\bf N}\left(\nabla\cdot{\bf N}\right) - \left({\bf N}\cdot\nabla\right){\bf N}\right]\cdot\nabla\Phi
\end{eqnarray}
where the additional term vanishes because ${\bf N}^2=1$ and ${\bf N}||\nabla\Phi$.  Recasting this in terms of the geometric invariants, we have
\begin{equation}
\vert\nabla\Phi\vert H = {1\over 2}\nabla\cdot\left\{\Phi\left[{\bf N}\left(\nabla\cdot{\bf N}\right) - \left({\bf N}\cdot\nabla\right){\bf N}\right]\right\} - K\Phi
\end{equation}
Thus, we can convert the cross term in (\ref{withcrossterm}) into a total derivative plus 
a purely geometric quantity~\cite{santangelo2}:
\begin{eqnarray}
F&=&{B\over 2}\int d^3\!x\,\Bigg[\left(1-\vert\nabla\Phi\vert - \lambda\nabla\cdot{\bf N}\right)^2 +4\lambda\Phi K
\nonumber\\
&&\quad+2\lambda\nabla\cdot\left\{{\bf N}-\Phi\left[{\bf N}\left(\nabla\cdot{\bf N}\right) - \left({\bf N}\cdot\nabla\right){\bf N}\right]\right\}\Bigg]
\end{eqnarray}
If we restrict ourselves to systems for which $K=0$, we again have the BPS decomposition.  Fortunately, edge dislocations introduce layer bending in only one direction, and hence a single
edge dislocation could be a BPS configuration, as could be a collection of parallel edge defects.  Whether or not such configurations exist depends on whether all boundary conditions of the problem can be fulfilled by the lower
order ``BPS'' equation:
\begin{equation}
1-\vert\nabla\Phi\vert = \lambda\nabla\cdot{\bf N}
\end{equation}

We have studied this equation numerically by switching to Lagrangian coordinates for the layers.  In terms of layers described by an embedding $\textbf{X}_n$ for layers labeled by $n$, this can be rewritten as the flow equation
\begin{equation}
{\bf N} \cdot \partial_n {\bf X}_n = \frac{a}{1 - 2 \lambda H}
\end{equation}
where we have used $\partial\Phi/\partial n =a$.
If we start with a layer at $n=0$ with vanishing Gaussian curvature, this flow generates new surfaces with vanishing Gaussian curvature, enabling us to find an exact solution for smectics with multiple parallel edge dislocations and to determine the energy due to bending the surrounding layers.  
\def\text#1{\hbox{#1}}
The nonlinear BPS equation precludes, on the surface, superposition.
However, in terms of $\tilde u =(z-\Phi)/\lambda$ and the new variables $\alpha=x/\sqrt{\lambda z}$ and
$\beta= z/\lambda$, we find an expansion around $\beta=\infty$~\cite{santangelo2}:
\begin{eqnarray}
0&=&\partial_\beta\tilde u- {1\over \beta}{1\over
2\left(1-\partial_\beta\tilde u\right)^3}\Bigg\{ \alpha\left(1-\partial_\beta\tilde u\right)^3\partial_\alpha\tilde u + \left[\left(1-\partial_\beta \tilde u\right)^2  -2\partial^2_\beta\tilde u\right]\left(\partial_\alpha\tilde u\right)^2\nonumber\\&&\quad\qquad -4\left(1-\partial_\beta\tilde u\right)\partial_\alpha\tilde u\partial_\alpha\partial_\beta\tilde u+2 \left(1-\partial_\beta\tilde u\right)^2\partial^2_\alpha\tilde u\Bigg\} +{\cal O}\left({1\over \beta^2}\right)
\end{eqnarray}
Thus, for large $\beta$, we may write our solutions as a power series around $\beta=\infty$.  The lowest order, non-trivial equation for $\tilde u=f(\alpha)$ is independent of $\beta$ and satisfies 
\begin{equation}\label{longdist}
{1\over 2}\left[\alpha\partial_\alpha f +\left(\partial_\alpha f\right)^2 \right] + \partial^2_\alpha f=0
\end{equation}
Through the Hopf-Cole transformation $S=\exp\{f/2\}$ we have
\begin{equation}\label{hc}
\alpha\partial_\alpha S +2 \partial^2_\alpha S=0,
\end{equation}
a linear equation.

Converting back to $y$ and $z$ coordinates, we find that $u$ approaches a similarity solution to the equation
\begin{equation}
\partial_z S - \lambda \partial_y^2 S = 0,
\end{equation}
for $z > 0$ and the same equation with opposite relative sign for $z < 0$, in analogy with the linear case.  The solution is~\cite{santangelo1,brener}
\begin{equation}\label{eq:BPSsoln}
u(y,z) = 2 \lambda \textrm{sgn}(z) \ln \left[1+\frac{e^{b/(4 \lambda)} - 1}{\sqrt{2 \pi}} \int_{-\infty}^{y/\sqrt{4 \lambda |z|}} e^{-t^2/2}~dt \right].
\end{equation}
The maximum strain lies along the parabolas defined by $y^2 = 4 \lambda |z|$.  Since, on these parabolas, $u(\sqrt{4 \lambda |z|},z)$ is constant for all $z$ it follows that equation~\ref{eq:BPSsoln} never approaches the linear solution along these parabolas even when $z \rightarrow \infty$.  In this sense, the linear theory never approaches the nonlinear theory.~\cite{brener,santangelo1}

\section{Nonzero Gaussian Curvature}

One difficulty of this approach is that it is no longer exact when $K \ne 0$.  The form of the BPS flow equation suggests, however, that we write the Euler-Lagrange equation for the extrema of the energy using $\Gamma = 1 - \vert\nabla \Phi \vert - \lambda \nabla \cdot {\bf N}$.  In this language, we find that~\cite{santangelo2}
\begin{equation}\label{eq:fullEL}
\nabla \cdot \left( {\bf N}~\Gamma - \lambda \frac{{\bf P}^T \cdot \nabla \Gamma}{\vert \nabla \Phi \vert} \right) = 2 \lambda K,
\end{equation}
where ${\bf P}^T={\hbox{$\mathbb{I}$}} - {\bf N}\otimes{\bf N}$ is the projection operator into the tangent plane of the layers.  When $K=0$, it is clear that $\Gamma=0$ is a solution to equation~\ref{eq:fullEL} as expected.  In addition, there will be $\Gamma \ne 0$ solutions corresponding to boundary conditions on $u$ that differ from those allowed by the BPS construction.

Equation~\ref{eq:fullEL} has the form of a conservation law for $\Gamma$, in which Gaussian curvature behaves as a source for $\Gamma$.  In general, this pair of equations is quite difficult to solve, and a full analysis of equation~\ref{eq:fullEL} has not yet been undertaken.  Nevertheless, there is one case with $K \ne 0$ that has been solved exactly in the nonlinear theory: the screw dislocation.  In the case of a screw dislocation, ${\bf m}= n a \hat z\delta(x)\delta(y)$ and so $\sigma(q)$ (from (\ref{tw1})) vanishes and 
\begin{equation}
u_{\rm screw} = -\frac{b}{2\pi}\tan^{-1}\left(\frac{y}{x}\right)
\end{equation}
in the quadratic theory.  However, taking $\Phi=z-u$ in the full theory, we find ($r^2=x^2+y^2$)
\begin{eqnarray}
{\bf N} &=& \frac{[-by/r,bx/r,2\pi r]}{\sqrt{b^2+(2\pi r)^2}}\nonumber\\
\Gamma &=&  1 - \sqrt{1+ b^2/(2\pi r)^2}
\end{eqnarray}
Note that $\nabla\cdot{\bf N}=0$ because the helicoid is a minimal surface.  Evaluating (\ref{eq:fullEL}), we see that the helicoid also is an extremum of the nonlinear theory.  Multiple screw dislocations are currently under investigation \cite{kamien,kamien3,santangelo3}.

\section{Summary and Open Problems}
The role of nonlinear elasticity in smectic liquid crystals is still a nascent field, and there are a number of open problems that are relatively unexplored.

Understanding the arrangement of layers in a particular smectic material is challenging precisely because of the interplay and frustration between the geometry of the layers and their spacing.  Though this problem can be solved exactly in some special cases, it is difficult to construct even reasonable conjectures for the shape of the layers in more general situations.  Nonetheless, the few known exact solutions available compel us to conclude that the harmonic approximation for the energy of smectic liquid crystals fails dramatically.

For the special case of $K=0$, exact solutions can be found that satisfy an equation of lower order than the Euler-Lagrange equation~\cite{santangelo1,santangelo2}.  These low order equations allow the determination of the shape of all the layers once the shape of a single layer has been specified.  Furthermore, the energy of these configurations is reduced to a topological term, which can be evaluated on the layers near the defect core.

Some analytical progress for $K \ne 0$ surfaces has also been made by directly summing the phase fields of screw dislocations.  In the case of a twist-grain boundary, where the screw dislocations lie along the line $x=0$, there is an unexpected connection to minimal surfaces~\cite{kamien,santangelo3}.  Neither sums of screw dislocations, nor their minimal surface counterparts, however, are exact minima of the smectic energy -- the true minimum likely lies somewhere in between~\cite{schick}.  More generally we can ask if there are solutions with $K\ne 0$ with energy determined entirely by the boundary conditions?

An analysis of the BPS flow for $\Gamma$ in the presence of Gaussian curvature has not yet been undertaken.  How do the layers evolve in the presence of Gaussian curvature?  Even when $K=0$, there are nontrivial aspects of the layer structure.  When there are two edge dislocations separated by several intervening layers, what is the structure of the layers in between and how do those edge dislocations interact?  In other words, can the BPS decomposition shown here be extended to solutions with constant or varying $K$?

In addition, because both screw and edge dislocations are topologically identical, there exists a family of dislocations with the same topology which can be found by rotating the defect core from being parallel (edge) to perpendicular (screw) to the layers.  These defects will have both edge and screw components, and we expect that nonlinear elasticity will be important in describing both their energetics (as is the case in screw dislocations) and the geometry of the smectic layers (as in edge dislocations).  Thus far, these defects have resisted theoretical descriptions beyond the harmonic approximation.  What is the structure and energetics of defects with both screw and edge components?

To summarize some open questions:
\begin{itemize}[$\bullet$]
\item At long distances the extremals of (\ref{functional}) satisfied (\ref{longdist}) which is actually a linear equation.  Thus, in this limit, we could find solutions for an arbitrary number of parallel edge dislocations 
at the {\sl same} $z$ co\"ordinate by setting the discontinuity between two solutions of the
diffusion equation, one diffusing ``up'' the other ``down'' \cite{santangelo1}.  >Can the BPS decomposition be used to
find solutions for parallel edge dislocations separated along the $z$ direction?

\item Single screw dislocations are helicoids and thus $H=0$, suggesting solutions
with vanishing compression, {\sl i.e.} $u_{zz}=0$.  However, $u_{zz}\ne 0$ for single
screw dislocations and so it is unlikely that there are superposition solutions which
satisfy $u_{zz}=0$.  We are free to add $F_K$ (\ref{gauss}) to (\ref{functional}):
\begin{equation}\label{test}
F=\frac{B}{2}\int d^3\!x\,\left\{\left(1-\vert\nabla\Phi\vert\right)^2 + \bar\lambda^2 K + \lambda^2 H^2\right\}
\end{equation}
with $\bar\lambda^2 = \bar K/B$.  >Is there a BPS decomposition of this functional for $H=0$?  >Is there a decomposition if we use a different power of $K$ in (\ref{test}), {\sl e.g.} $K^2$?

\item >Are there extremals of (\ref{functional}) which are not parallel to the $z$ axis (screw dislocations) or lying in the $xy$-plane (edge dislocations)?  >How does the BPS edge solution deform into the minimal surface screw dislocation?
\end{itemize}

The solution of these problems has the added benefit that they can be verified directly through experiments on liquid crystals~\cite{oleg} and the physical validity of any approximations can
thus be assessed.  Despite the long history of theory and experiment on these layered systems, the nonlinear aspects of these problems are only now being explored. 

\acknowledgements
We have benefitted from discussions with Gregory Grason and Tom Lubensky.  This work was supported in part by NSF Grants DMR01-29804 and DMR05-47230, the Donors of the Petroleum Research Fund, and a gift from Lawrence J. Bernstein.

% The endnotes section will be placed here.

%\theendnotes

\end{article}

\begin{thebibliography}{}

\bibitem[\protect\citeauthoryear{Bogomol'nyi}{1976}]{BPS1}
E.B. Bogomol'nyi
\newblock {Stability of Classical Solutions}.
\newblock {\em Sov. J. Nucl. Phys.}, 24:449--454, 1976.

\bibitem[\protect\citeauthoryear{Brener and Marchenko}{1999}]{brener}
E.A. Brener and V.I. Marchenko
\newblock {Nonlinear theory of dislocations in smectic crystals: an exact solution}.
\newblock {\em Phys. Rev. E}, 59:R4752, 1999.

\bibitem[\protect\citeauthoryear{Chaikin and Lubensky}{1995}]{chaikin}
P.M. Chaikin and T.C. Lubensky
\newblock {\textit{Principles of Condensed Matter Physics}}.
\newblock (Cambridge University Press, New York, 1995).

\bibitem[\protect\citeauthoryear{de Gennes}{1972}]{degennes72}
P.G. de Gennes
\newblock{Analogy between superconductors and smectics-\textit{A}}.
\newblock {\em Solid State Commun.}, 10:753--756, 1972.

\bibitem[\protect\citeauthoryear{de Gennes and Prost}{1993}]{degennes}
P.G. de Gennes and J. Prost
\newblock {\textit{The Physics of Liquid Crystals}}.
\newblock (Oxford University Press, New York, 1993).

\bibitem[\protect\citeauthoryear{DiDonna and Kamien}{2003}]{didonna}
B.A. DiDonna and R.D. Kamien
\newblock {Smectic Blue Phases: Layered Systems with High Intrinsic Curvature}.
\newblock {\em Phys. Rev. E}, 68:041703, 2003.

\bibitem[\protect\citeauthoryear{Duque and Schick}{2000}]{schick}
D. Duque and M. Shick
\newblock {Self-consistent field theory of twist grain boundaries in block copolymers}.
\newblock {\em J. Chem. Phys.}, 113:5525--5530, 2000.

\bibitem[\protect\citeauthoryear{Ishikawa and Lavrentovich}{1999}]{oleg}
T. Ishikawa and O.D. Lavrentovich
\newblock {Dislocation Profile in Cholesteric Finger Texture}.
\newblock {\em Phys. Rev. E}, 60:R5037Ð5039, 1999.

\bibitem[\protect\citeauthoryear{Kamien and Lubensky}{1999}]{kamien}
R.D. Kamien and T.C. Lubensky
\newblock {Minimal Surfaces, Screw Dislocations, and Twist-Grain Boundaries}.
\newblock {\em Phys. Rev. Lett.}, 82:2892--2895, 1999.

\bibitem[\protect\citeauthoryear{Kamien}{1999}]{kamien3}
R.D. Kamien
\newblock {Decomposition of the Height Function of Scherk's First Surface}.
\newblock {\em Appl. Math. Lett.}, 14:797--800, 2001.

\bibitem[\protect\citeauthoryear{Kamien}{2002}]{kamien2}
R.D. Kamien
\newblock {The Geometry of Soft Materials: A Primer}.
\newblock {\em Rev. Mod. Phys.}, 74:953--972, 2002.

\bibitem[\protect\citeauthoryear{Kleman}{1983}]{klemanbook}
M. Kleman
\newblock {\textit{Points, Lines, and Walls}}.
\newblock (John Wiley Sons, New York, 1983).

\bibitem[\protect\citeauthoryear{Prasad and Sommerfield}{1975}]{BPS2}
M.K. Prasad and C.H. Sommerfield
\newblock {Exact Classical Solution for the 't Hooft Monopole and the Julia-Zee Dyon}.
\newblock {\em Phys. Rev. Lett.}, 35:760-762, 1975.


\bibitem[\protect\citeauthoryear{Santangelo and Kamien}{2003}]{santangelo1}
C.D. Santangelo and R.D. Kamien
\newblock {Bogomol'nyi-Prasad-Sommerfield Configurations in Smectics}.
\newblock {\em Phys. Rev. Lett.}, 91:045506, 2003.

\bibitem[\protect\citeauthoryear{Santangelo and Kamien}{2005a}]{santangelo2}
C.D. Santangelo and R.D. Kamien
\newblock {Curvature and Topology in Smectic Liquid Crystals}.
\newblock {\em Proc. R. Roy. Soc. A}, 461:2911--2921, 2005.

\bibitem[\protect\citeauthoryear{Santangelo and Kamien}{2005b}]{santangelo3}
C.D. Santangelo and R.D. Kamien
\newblock {Elliptic Phases:  A Study of the Nonlinear Elasticity of Twist-Grain Boundaries}.
\newblock {\em cond-mat/0511740}, 2005.

\end{thebibliography}
\end{document}